\documentclass[11pt]{article}
\usepackage{amsmath,amssymb,amsthm,geometry,verbatim,enumerate,float,tikz}
\newtheorem{thm}{Theorem}
\newtheorem{lemma}[thm]{Lemma}

\newtheorem{claim}{Claim}

\usepackage{setspace}

\begin{document}

\onehalfspace

\title{Longest Paths in Circular Arc Graphs}

\author{Felix Joos}

\date{}

\maketitle

\begin{center}
Institut f\"{u}r Optimierung und Operations Research, 
Universit\"{a}t Ulm, Ulm, Germany\\
\texttt{felix.joos@uni-ulm.de}
\end{center}

\begin{abstract}
As observed by Rautenbach and Sereni (arXiv:1302.5503)
there is a gap in the proof of the theorem of Balister et al.
(Longest paths in circular arc graphs,
Combin. Probab. Comput., 13, No. 3, 311-317 (2004)),
which states that the intersection of all longest paths in a connected circular arc graph is nonempty.
In this paper we close this gap.

\bigskip

\end{abstract}

\section{Introduction}
It is easy to prove that every two longest paths in a connected graph have a nonempty intersection.
Gallai \cite{ga} asked if the intersection of all longest paths is nonempty.
This is not true in general but holds for some graph classes.
See \cite{zam} for a survey.
In \cite{ba} Balister et al. proved that it is true for interval graphs and circular arc graphs.
	However, as pointed out by Rautenbach and Sereni \cite{rs},
there is a gap in the proof for the class of circular arc graphs.
Rautenbach and Sereni proved the weaker result
that in a connected circular arc graph, there is a set of at most $3$ vertices such that every longest paths intersects this set.
We close the gap by extending Lemma 3.2 from \cite{ba}.
Our Lemma \ref{mainresult} corresponds to Lemma 3.2 in \cite{ba}.

We follow the notation in \cite{ba}.
A graph $G$ is a \textit{circular arc graph},
if there exists a function $\phi$ of its vertex set $V(G)$ into a collection of open arcs of a circle such that,
for every two distinct vertices $u$ and $w$ of $G$, $uw$ is an edge of $G$ if and only if $\phi(u)\cap \phi(w)\neq \emptyset$,
that is, the class of circular arc graphs are the intersection graphs of arcs in a circle.
Let \textit{interval graphs} be the intersection graphs of open intervals of the real line.
Note that one can assume that all endpoints of the arcs and intervals are distinct.

\section{Result}
We review the approach of Balister et al.
Let $G$ be a connected circular arc graph.
Let $C$ be a circle and $\mathcal{F}$ be a finite collection of open arcs of $C$
that correspond to the vertices of $G$.
If the union of arcs in $\mathcal{F}$ does not cover $C$,
then $G$ is an interval graph and hence the statement follows by a result of \cite{ba}.
Therefore, we may assume that the union of arcs in $\mathcal{F}$ covers $C$.
We choose a set $\mathcal{K}\subseteq \mathcal{F}$ such that
$\mathcal{K}=\{K_0,\ldots,K_{n-1}\}$,
\begin{itemize}
	\item $C=K_0\cup\ldots\cup K_{n-1}$,
	\item $n$ is minimal, and
	\item no $K_i$ is contained in another arc, i.e. $K_i\subseteq A\in \mathcal{F}\Rightarrow K_i=A$.
\end{itemize}
We cyclically order the elements of $\mathcal{K}$ clockwise and
consider all indices of elements of $\mathcal{K}$ modulo $n$.
A \textit{chain} $\mathcal{P}$ of length $t$ is a $t$-tupel $(J_1,\ldots,J_t)$ of distinct arcs (in $\mathcal{F}$) such that
$J_i\cap J_{i+1}\neq \emptyset$ for every $1\leq i \leq t-1$.
This corresponds to a path in $G$ on $t$ vertices.
The chain $\mathcal{P}$ is a \textit{longest} chain,
if there is no chain of larger length than $\mathcal{P}$.
For a chain $\mathcal{P}=(J_1,\ldots,J_t)$, 
let the \textit{support} ${\rm Supp\ }\mathcal{P}$ of $\mathcal{P}$ be the subset of $C$ defined by

$$J_1\cup(J_2\cap J_3)\cup \ldots \cup(J_{t-2}\cap J_{t-1})\cup J_t.$$
Note that if there is an arc $A$ in $\mathcal{F}$ that is not contained in the chain $\mathcal{P}$ of length $t$ and intersects ${\rm Supp\ }\mathcal{P}$,
then there is a chain of length $t+1$ consisting of the arc $A$ and all arcs of $P$.
This implies that for a longest chain $\mathcal{P}$ in $\mathcal{F}$, an arc $A$ is contained in $\mathcal{P}$
if and only if it intersects ${\rm Supp\ } \mathcal{P}$.

For two points $x,y$ on the circle $C$, let $[x,y]$ be the arc from $x$ to $y$ in clockwise direction.
For an arc $A\in \mathcal{F}$,
let $\ell(A)$ and $r(A)$ be the left and right endpoint of $A$, respectively,
that is $\ell(A),A,r(A)$ are consecutive on $C$ in clockwise direction.

Now, we mention two results, which we use later.

\begin{lemma}[Balister et al. \cite{ba}]\label{consec}
If $\mathcal{P}$ is a longest chain in $\mathcal{F}$,
then $\mathcal{P}\cap \mathcal{K}=\{K_i:i\in I\}$ is nonempty
and $I$ is a contiguous set of elements of $\mathbb{Z}_n$.
\end{lemma}

\noindent
The next lemma is due to Keil \cite{ke} and explicitly formulated as Lemma 2.3 in \cite{ba}.

\begin{lemma}[Keil \cite{ke}]\label{reorder}
Let $X=\{x_1,\ldots,x_{t+1}\}$ be a set of real numbers, 
and let $J_1,\ldots,J_t$ be a sequence of open intervals with $x_k,x_{k+1}\in J_k$ for every $1\leq k \leq t$.
If $x_{i_1}<\ldots<x_{i_{t+1}}$ are the elements of $X$ in increasing order,
then the intervals have a permutation $J_{j_1},\ldots,J_{j_t}$ such that $x_{i_k},x_{i_{k+1}}\in J_{j_k}$, for every $1\leq k \leq t$.
\end{lemma}

\begin{figure}[t]
\begin{center}
\begin{tikzpicture}[scale=1]
\def\ver{0.1} 

\draw[thick] (3,0)--(9,0)
(1,1)--(8,1)
(0,2)--(6,2)
(-1,3)--(2,3);

\draw (7,1.3) node {$J_q$};
\draw (8,0.3) node {$J_{q+1}$};
\draw (0.5,2.3) node {$J_p$};
\draw (0,3.3) node {$J_{p-1}$};

\fill (1.5,3) circle (\ver);
\fill (1.5,2) circle (\ver);
\fill (1.5,1) circle (\ver);
\draw (1.5,3.3) node {$x_{p}$};

\fill (2.5,2) circle (\ver);
\fill (2.5,1) circle (\ver);
\draw (2.5,2.3) node {$x_{p+1}$};

\fill (5.5,2) circle (\ver);
\fill (5.5,1) circle (\ver);
\fill (5.5,0) circle (\ver);
\draw (5.5,2.3) node {$x_{q+1}$};

\fill (4.5,2) circle (\ver);
\fill (4.5,1) circle (\ver);
\draw (4.5,2.3) node {$x_{q}$};

\end{tikzpicture}
\end{center}
\caption{$J_p$ and $J_q$ can be swapped.}
\label{swap}
\end{figure}
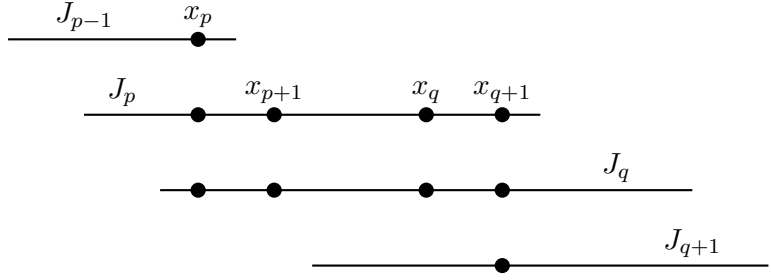

Let $\mathcal{P}=(J_1,\ldots, J_t)$ be a chain such that $\mathcal{K}\not\subseteq \mathcal{P}$
and let $\{x_1,\ldots,x_{t+1}\} \subset$ ${\rm Supp\ } \mathcal{P}$ be a set of distinct points such that $x_k,x_{k+1}\in J_k$,
for every $1\leq k\leq t$.
Without loss of generality, we may assume, by Lemma \ref{reorder}, that $x_1,x_2,\ldots,x_{t+1}$ are consecutive points on $C$ in clockwise direction.
One might have to replace $\mathcal{P}$ by another chain having exactly the same arcs.
Let $p,q\in \{1,\ldots,t\}$ such that $p<q$.
If $[x_p,x_{p+1}],[x_q,x_{q+1}]\subseteq J_p\cap J_q$,
then the reordering $$(J_1,\ldots, J_{p-1}, J_q, J_{p+1},\ldots, J_{q-1}, J_p, J_{q+1},\ldots,J_t)$$ of $\mathcal{P}$ is a chain of the same length as $\mathcal{P}$.
See Figure \ref{swap} for illustration.
In this situation it is possible to \textit{swap} $J_p$ and $J_q$ in $\mathcal{P}$.

For $i\in \{0,\ldots,n-1\}$,
let $\Delta K_i=\{x\in C: \ell(K_{i+1})<x<r(K_{i})\}$.
Note that for $n\geq 3$,
we have $\Delta K_i=K_i \cap K_{i+1}$.
We use this notation because Balister et al. omitted the case $n=2$.
Note that $A\subset K_{i}\cup K_{i+1}$ implies the connectedness of $A\setminus \Delta K_{i+1}$ if $n$ is at least 3.

Lemma \ref{mainresult} is our main contribution.
Balister et al. only proved Lemma \ref{mainresult} with the properties $(a)$ - $(c)$.
We extend this result.

\begin{lemma}\label{mainresult}
If $\mathcal{P}$ is a longest chain in $\mathcal{F}$ and $\mathcal{P}\cap \mathcal{K}=\{K_{a+1},\ldots,K_{b-1}\}\neq \mathcal{K}$,
then the arcs in $\mathcal{P}$ have a reordering into a chain $\mathcal{P}^*$ such that in this reordering
\begin{enumerate}[(a)]
	\item $K_{a+1}$ precedes $K_{b-1}$ in $\mathcal{P}^*$ provided they are distinct,
	\item if $A$ precedes $K_{b-1}$ in $\mathcal{P}^*$, then $\Delta K_{b-1} \not\subseteq A$,
	\item if $A$ precedes $K_{a+1}$ in $\mathcal{P}^*$, then $A\subseteq K_a\cup K_{a+1}$ and $A\setminus \Delta K_{a+1}$ is connected,
	\item if $K_{b-1}$ precedes $A$ in $\mathcal{P}^*$, then $A\subseteq K_{b-1}\cup K_{b}$ and $A\setminus \Delta K_{b}$ is connected,
	\item if $K_{a+1}$ precedes $A$ in $\mathcal{P}^*$, then $\Delta K_{a} \not\subseteq A$.
\end{enumerate}
\end{lemma}

Here is the gap of Balister et al. 
Indeed (b) and (c) is symmetric to (d) and (e) (they proved that (b) and (c) holds), 
however,
forcing both at the same time is a stronger assertion.

\bigskip

\noindent
\textit{Proof:}
Let $\mathcal{P}=(J_1,\ldots, J_t)$ 
and let $\{x_1,\ldots,x_{t+1}\} \subset$ Supp $\mathcal{P}$ be a set of distinct points such that $x_k,x_{k+1}\in J_k$
for every $1\leq k\leq t$.
Without loss of generality, we may assume, by Lemma \ref{reorder}, that $x_1,x_2,\ldots,x_{t+1}$ are consecutive points on $C$ in clockwise direction.
It is important to keep in mind that every $x_i$ belongs to $(K_{a+1}\cup \ldots \cup K_{b-1})\setminus (K_a\cup K_b)$,
because $K_a$ and $K_b$ do not belong to $\mathcal{P}$.

First, we prove $(c)$ and $(e)$.
Let $\mathcal{P}'=(J_{j_1},\ldots,J_{j_s})$ be a subsequence of $\mathcal{P}$ such that $A\in \mathcal{P}'$
if and only if
\begin{enumerate}[(i)]
	\item $K_{a+1}$ precedes $A$ and $\Delta K_{a} \subseteq A$ or 
	\item $A$ precedes $K_{a+1}$ and $A\not\subseteq K_a \cup K_{a+1}$ if $n\geq3$ and $A\setminus \Delta K_{a+1}$ is disconnected if $n=2$.
\end{enumerate}

If $n\geq 3$, then we observe the following.
If $A\in\mathcal{P}'$ satisfies requirement (i),
then, by the choice of $\mathcal{K}$, we conclude that $\ell(K_a)$, $\ell(A)$, $\ell(K_{a+1})$, $r(K_a)$, $r(A)$, $r(K_{a+1})$ are consecutive points in clockwise direction on $C$.
If $A\in\mathcal{P}'$ satisfies requirement (ii),
then $\ell(K_a)$, $\ell(K_{a+1})$, $\ell(A)$, $r(K_a)$, $r(K_{a+1})$, $r(A)$ 
or $\ell(K_a)$, $\ell(K_{a+1})$, $r(K_a)$, $\ell(A)$, $r(K_{a+1})$, $r(A)$ are consecutive points in clockwise direction on $C$,
because $A\cap (K_{a+1}\setminus K_a)\not=\emptyset$ and the choice of $\mathcal{K}$.
For $n=2$ the situation is obvious.
Let 
\begin{itemize}
	\item $L=\{i\in[t]: J_i\in \mathcal{P} {\rm\ and\ }J_i{\rm\ satisfies\ requirement\ (i)}\}$ and
	\item $R= \{i\in[t]: J_i\in \mathcal{P} {\rm\ and\ }J_i{\rm\ satisfies\ requirement\ (ii)}\}$.
\end{itemize}
Let $L_{\mathcal{P}}=\{J_i \in \mathcal{P}:i\in {\rm L}\}$ and $R_{\mathcal{P}}=\{J_i \in \mathcal{P}:i\in {\rm R}\}$,
that is $L_{\mathcal{P}}$ and $R_{\mathcal{P}}$ partition $\mathcal{P}'$.
Furthermore, all arcs in $R_{\mathcal{P}}$ precede the arcs in $L_{\mathcal{P}}$.
Note that all arcs in $\mathcal{P}\setminus \mathcal{P}'$ satisfy the requirements (c) and (e).
Let $\gamma\in\mathbb{N}$ be such that $K_{a+1}=J_{\gamma}$ and  $f(\mathcal{P}')$ be defined by 
$$\max\{\{\gamma\}\cup L \cup R \}-\min\{\{\gamma\}\cup L \cup R \}.$$

\begin{claim}\label{c1}
Let $L$ and $R$ be non empty, and consider $p\in R$ and $q\in L$.
Then it is possible to swap $J_p$ and $J_q$ in $\mathcal{P}$, the reordering of $\mathcal{P}$ is still a chain
and the sets $L$ and $R$ lose exactly $q$ and $p$, respectively.
\end{claim}

\noindent
\textit{Proof of Claim \ref{c1}:}
By our observations above and since $J_p$ precedes $J_q$,
we conclude that $\ell(J_q)$, $\ell(J_p)$, $r(J_q)$ and $r(J_p)$ are consecutive points in clockwise direction on $C$.
Since $J_p$ precedes $J_q$, 
we obtain $[x_{p},x_{p+1}],[x_{q},x_{q+1}]\subseteq J_p \cap J_q$.
Thus it is possible to swap $J_p$ and $J_q$ in $\mathcal{P}$.
After this swap both arcs do not satisfy the requirements (i) and (ii) any more
and in addition the relative positions of all other arcs concerning $K_{a+1}$ do not change.
This completes the proof of the claim.
$\Box$

\begin{claim}\label{c2}
Each element $J_p\in R_{\mathcal{P}}$ can be swapped with $K_{a+1}$ in $\mathcal{P}$ and 
the reordering of $\mathcal{P}$ is still a chain.
\end{claim}

\noindent
\textit{Proof of Claim \ref{c2}:}
Let $q$ be such that $J_q=K_{a+1}$, that is $p<q$ by the definition of $R$.
By our observations above,
we know that $\ell(K_{a+1})$, $\ell(J_p)$, $r(K_{a+1})$ and $r(J_p)$ are consecutive points in clockwise direction on $C$.
Since $J_p$ precedes $K_{a+1}$,
we obtain $[x_p,x_{p+1}],[x_q,x_{q+1}]\subseteq J_p\cap K_{a+1}$.
Thus it is possible to swap $J_p$ and $K_{a+1}$ in $\mathcal{P}$ and the reordering of $\mathcal{P}$ is still a chain.
$\Box$

\begin{claim}\label{c3}
Each element $J_q\in L_{\mathcal{P}}$ can be swapped with $K_{a+1}$ in $\mathcal{P}$ and 
the reordering of $\mathcal{P}$ is still a chain.
\end{claim}

\noindent
\textit{Proof of Claim \ref{c3}:}
Let $p$ be such that $J_p=K_{a+1}$, that is $p<q$ by the definition of $L$.
By our observations above,
we know that $\ell(J_q)$, $\ell(K_{a+1})$, $r(J_q)$ and $r(K_{a+1})$ are consecutive points on in clockwise direction $C$.
Since $K_{a+1}$ precedes $J_q$,
we obtain $[x_p,x_{p+1}],[x_q,x_{q+1}]\subseteq K_{a+1} \cap J_q$.
Thus it is possible to swap $K_{a+1}$ and $J_q$ and the reordering of $\mathcal{P}$ is still a chain.
$\Box$

\bigskip

Recall that $J_{\gamma}=K_{a+1}$.
Let $\alpha=\min\{\{\gamma\}\cup L \cup R \}$ and $\beta=\max\{\{\gamma\}\cup L \cup R \}$.
Note that $\alpha$ does not decrease and $\beta$ does not increase
if we reorder $\mathcal{P}$ as described in Claims \ref{c1}-\ref{c3}.
In particular, $f(\mathcal{P}')$ does not increase.
After swapping two elements in $\mathcal{P}'$, by Claim \ref{c1} the subsequence loses two elements.
Using Claim \ref{c1} iteratively, we can assume that $L=\emptyset$ or $R=\emptyset$.
If $\mathcal{P}'=\emptyset$, then this completes the proof of (c) and (e).
Therefore, we assume that $\mathcal{P}'\not=\emptyset$ and $\mathcal{P}'=L_{\mathcal{P}}$ or $\mathcal{P}'=R_{\mathcal{P}}$.
We distinguish the two possible cases.
\begin{enumerate}[I]
	\item If $\mathcal{P}'=L_{\mathcal{P}}$, 
				then we have $K_{a+1}=J_\alpha$ and $\beta=\max\{L\}$, and
	\item if $\mathcal{P}'=R_{\mathcal{P}}$,
				then we have $\alpha=\min\{R\}$ and $K_{a+1}=J_\beta$.
\end{enumerate}
Note that $f(\mathcal{P}')=0$ if and only if $\mathcal{P}'=\emptyset$.
By Claims 2 and 3, it is possible to swap $K_{a+1}$ with each element of $\mathcal{P}'$.
In the first case swap $K_{a+1}$ with $J_{\beta}$ and in the second case with $J_{\alpha}$.
Denote this reordering of $\mathcal{P}$ by $\mathcal{P}$ again and define $\mathcal{P}'$, $L$ and $R$ as before.
Consider first case I.
Note that $L=\emptyset$ and $R\subseteq \{\alpha+1,\ldots,\beta-1\}$.
In case II, we have $L\subseteq \{\alpha+1,\ldots,\beta-1\}$ and $R=\emptyset$.
In both cases $f(\mathcal{P}')$ decreases by at least 1.
After iterating this procedure at most $\beta-\alpha$ times,
we have $f(\mathcal{P}')=0$.
Hence there is a reordering of $\mathcal{P}$ such that the requirements $(c)$ and $(e)$ are fulfilled.
From now on, we assume that $\mathcal{P}$ fulfills requirements $(c)$ and $(e)$.

If $a+1=b-1$,
then $\mathcal{P}$ fulfills the requirements (a), (b) and (d).
Note that this is also true if $|\mathcal{K}|=2$.
Thus we assume that $K_{a+1}$ and $K_{b-1}$ are distinct.
This implies $n\geq 3$.
Note that $K_{a+1}$ precedes $K_{b-1}$ by requirement (c).
Let $\tilde{\mathcal{P}}=(J_{k_1},\ldots,J_{k_{s'}})$ be the subsequence of $\mathcal{P}$ such that $A\in \tilde{\mathcal{P}}$
if and only if
\begin{enumerate}[(i')]
	\item $K_{b-1}$ precedes $A$ and $A\not\subseteq \Delta K_{b-1}$ or 
	\item $A$ precedes $K_{b-1}$ and $K_{b-1}\cap K_{b} \subseteq A$.		
\end{enumerate}
Note that $K_{a+1}\notin \tilde{\mathcal{P}}$.
Let $\tilde{L}=\{i\in[t]: J_i\in \mathcal{P} {\rm\ and\ }J_i{\rm\  satisfies\ requirement\ (i')}\}$ and
$\tilde{R}= \{i\in[t]: J_i\in \mathcal{P} {\rm\ and\ }J_i{\rm\  satisfies\ requirement\ (ii')}\}$.
Let $\tilde{\gamma}\in\mathbb{N}$ be such that $K_{b-1}=J_{\tilde{\gamma}}$ and
$\tilde{\alpha}=\min\{\{\tilde{\gamma}\}\cup \tilde{L} \cup \tilde{R} \}$.
Note that $\gamma<\tilde{\alpha}$.
This implies that $K_{a+1}$ precedes all arcs in $\tilde{\mathcal{P}}$ and
hence arguing as above for $K_{b-1}$, the relative order in the ordering of $\mathcal{P}$ of all arcs of $\mathcal{P}$ 
concerning $K_{a+1}$ does not change.
This shows that there is a reordering $\mathcal{P}^*$ of $\mathcal{P}$ such that $\mathcal{P}^*$ fulfills the requirements of Lemma \ref{mainresult}.
$\Box$.

\begin{thm}\label{mainthm}
If $G$ is a connected circular arc graph,
then the intersection of all longest paths is nonempty.
\end{thm}

\noindent
\textit{Proof:} We can assume that $G$ is not an interval graph, otherwise the statement follows by a result of \cite{ba}.
As above, let $\mathcal{F}$ be the finite collection of arcs of a circle $C$
that correspond to the vertices of $G$.
We choose $\mathcal{K}$ as above.
If $n=1$, then every longest chain contains $K_0$ and we are done.
Let $\mathcal{P}$ a longest chain such that $|\mathcal{P}\cap \mathcal{K}|$ is as small as possible.
If $|\mathcal{P}\cap \mathcal{K}|=n$,
then every longest chain contains all arcs of $\mathcal{K}$ and we are done, too.
Therefore, we assume that $n\geq 2$ and $|\mathcal{P}\cap \mathcal{K}|<n$.
That is, by Lemma \ref{consec}, 
$\mathcal{P}\cap \mathcal{K}=\{K_{a+1},\ldots,K_{b-1}\}$.
We prove Theorem \ref{mainthm} by showing that every longest chain contains $K_{b-1}$.
We assume for contradiction, 
that there is a longest chain $\mathcal{Q}$ such that $K_{b-1}\notin \mathcal{Q}$.
Let $\mathcal{Q}\cap \mathcal{K}=\{K_{\ell+1},\ldots,K_{m-1}\}$.
Our assumption and choice of $\mathcal{P}$ imply that 
$K_{b-1}\in \mathcal{P}\setminus \mathcal{Q}$, $K_{\ell+1}\in \mathcal{Q}\setminus \mathcal{P}$ and
$K_{b},\ldots,K_{\ell}\notin \mathcal{P}\cup\mathcal{Q}$.
Let $\mathcal{R}$ be the chain $(K_b,\ldots,K_{\ell}$).
Note that $\mathcal{R}=\emptyset$, if $b=\ell+1$.

For a $k$-tuple $\mathcal{A}=(A_1,\ldots,A_k)$,
let the reversed $k$-tuple $\mathcal{A}^r$ be defined by $(A_k,\ldots,A_1)$.
If $\mathcal{B}=(B_1,\ldots,B_{k'})$,
then let $\mathcal{A}\mathcal{B}=(A_1,\ldots,A_k,B_1,\ldots,B_{k'})$ and $\mathcal{A}B_1=(A_1,\ldots,A_k,B_1)$.
We reorder $\mathcal{P}$ and $\mathcal{Q}$ such that the reorderings $\mathcal{P}^*$ and $\mathcal{Q}^*$ satisfy the requirements of Lemma \ref{mainresult}.
Let $\mathcal{P}^*=\mathcal{P}_1 K_{b-1} \mathcal{P}_2$ and $\mathcal{Q}^*=\mathcal{Q}_1 K_{\ell+1} \mathcal{Q}_2$.
Note that 
\begin{enumerate}[(i)]
	\item if $A\in \mathcal{P}_1$, then $\Delta K_{b-1} \not\subseteq A$,
	\item if $A\in \mathcal{P}_2$, then $A\subseteq K_{b-1}\cup K_{b}$ and $A\setminus \Delta K_{b}$ is connected,
	\item if $A\in \mathcal{Q}_1$, then $A\subseteq K_{\ell}\cup K_{\ell+1}$ and $A\setminus \Delta K_{\ell+1}$ is connected, and
	\item if $A\in \mathcal{Q}_2$, then $\Delta K_{\ell} \not\subseteq A$.
\end{enumerate}
Let $\mathcal{C}_1=\mathcal{P}_1K_{b-1}\mathcal{R}K_{\ell+1}\mathcal{Q}_1^r$ and
$\mathcal{C}_2=\mathcal{P}_2^rK_{b-1}\mathcal{R}K_{\ell+1}\mathcal{Q}_2$.	

\setcounter{claim}{0}

\begin{claim}
$\mathcal{C}_1$ is a chain.
\end{claim}

\noindent
\textit{Proof:} It suffices to show that $\mathcal{P}_1\cap\mathcal{Q}_1=\emptyset$.
We assume for contradiction, that there is an arc $A\in \mathcal{P}_1\cap \mathcal{Q}_1$.
Suppose $n=2$.
Thus $\mathcal{K}=\{K_{b-1},K_{\ell+1}\}$.
By (iii), $A\setminus \Delta K_{\ell+1}$ is connected and by (i) $\Delta K_{b-1} \not\subseteq A$.
This implies that $A\subseteq K_{\ell+1}$ or $A\subseteq K_{b-1}$.
Since $A\in\mathcal{P}\cap \mathcal{Q}$,
this implies $K_{\ell+1}\in \mathcal{P}\cap \mathcal{Q}$
or $K_{b-1}\in \mathcal{P}\cap \mathcal{Q}$,
which is a contradiction.

Now we assume $n\geq 3$.
By (iii), $A\subseteq K_{\ell}\cup K_{\ell+1}$.
Since $A$ meets $K_{\ell+1}\setminus K_{\ell}$,
we observe that $r(K_{\ell})$, $r(A)$ and $r(K_{\ell+1})$ are consecutive points on $C$.
If $A\subseteq K_{\ell+1}$,
then ${\rm Supp\ }\mathcal{P}\cap K_{\ell+1}\not=\emptyset$ and hence $K_{\ell+1}\in \mathcal{P}$, which is a contradiction.
Thus $\ell(K_{\ell})$, $\ell(A)$, $\ell(K_{\ell+1})$, $r(K_{\ell})$, $r(A)$ and $r(K_{\ell+1})$ are consecutive points on $C$.

\noindent
By (i), $K_{b-1}\cap K_b \not\subseteq A$.
This implies that $b\not=\ell+1$ and hence $\mathcal{R}$ is not empty.
Thus $K_{\ell}\notin \mathcal{P}$.
Since $A\in \mathcal{P}$, 
it is $A\cap{\rm Supp\ }\mathcal{P}\not= \emptyset$ and hence ${\rm Supp\ }\mathcal{P}\cap (K_{\ell}\cup K_{\ell+1})\not= \emptyset$.
Thus $\mathcal{P}$ contains $K_{\ell}$ or $K_{\ell+1}$.
This is a contradiction and completes the proof of Claim 1.
$\Box$

\begin{claim}
$\mathcal{C}_2$ is a chain.
\end{claim}

\noindent
\textit{Proof:} Using (ii) and (iv) instead of (i) and (iii) this is the completely symmetric case to Claim 1.
$\Box$

\bigskip

\noindent
Note that $|\mathcal{C}_1|+|\mathcal{C}_2|\geq |\mathcal{P}|+|\mathcal{Q}|+2$.
This implies that $|\mathcal{C}_1|>|\mathcal{P}|$ or $|\mathcal{C}_2|>|\mathcal{P}|$,
which is a contradiction to the choice of $\mathcal{P}$.
$\Box$

\end{document}